\documentstyle[12pt]{article}

  \newcommand{\Cov}{\rm  Cov}

\begin{document}

   \begin{center}

{\bf Central Limit Theorem in Lebesgue-Riesz spaces  }\\

\vspace{4mm}

{\bf  for weakly dependent random sequences}

 \vspace{5mm}

 {\bf M.R.Formica, \ E.Ostrovsky, \ L.Sirota }

 \end{center}

\vspace{5mm}

 \ Universit\`{a} degli Studi di Napoli Parthenope, via Generale Parisi 13, Palazzo Pacanowsky, 80132,
Napoli, Italy \\

e-mail: mara.formica@uniparthenope.it \\

\vspace{4mm}

 \ Israel,  Bar-Ilan University, department of Mathematic and Statistics, 59200 \\

\vspace{4mm}

   e-mail:eugostrovsky@list.ru\\
   
   e-mail:sirota3@bezeqint.net \\

\vspace{5mm}

\begin{center}

 {\bf Abstract} \\

\end{center}

\vspace{5mm}

 \ We deduce sufficient conditions for the Central Limit Theorem (CLT) in the Lebesgue-Riesz space $ \ L(p) \ $
 defined on some measure space for the sequence of centered random variables satisfying the strong mixing (Rosenblatt) condition. \par

  \ We investigate the asymptotical as well as non-asymptotical approach. \par

\vspace{5mm}

\begin{center}

 {\bf Key words and phrases} \\

\end{center}

\vspace{5mm}

 \ Random variables and vectors (r.v.), distributions, sigma algebras (fields), probability, estimations, metrisable spaces, separability, Utev's constant,
 compactness, normed sums, Central Limit Theorem (CLT), Lebesgue-Riesz space, (strong) mixing conditions and
 mixing (Rosenblatt) coefficient, weak convergence of distributions.

\vspace{5mm}

\section{Statement of problem, definitions, previous results.  }

\vspace{5mm}

 \ Let $ \ (T = \{t\},B,\mu) \ $ be measurable space with sigma-finite separable non-zero measure $ \ \mu. \ $ The separability implies that the  sigma-field $ \ B \ $ is separable
 metric space relative the distance function

 $$
 d(A_1,A_2) := \mu(A_1 \setminus A_2) + \mu(A_2 \setminus A_1).
 $$

 \ Recall that the classical Lebesgue-Riesz space $ \ L(p) = L(p,T,\mu) = L(p,T) \ $  consists on all the measurable numerical valued functions $ \ f: T \to R \ $ having a finite norm

$$
||f||_{p,T}  = ||f||L(p,T) = ||f||L(p,T,\mu) \stackrel{def}{=} \left[ \ \int_T |f(t)|^p \ \mu(dt) \  \right]^{1/p}, \ 1 \le p < \infty.
$$

 \ It is a complete separable rearrangement invariant (r.i.) Banach functional space. \par

 \ Let also $ \ (\Omega = \{\omega \}, F,P) \ $ be non-trivial probability space. The Lebesgue-Riesz norm of arbitrary numerical valued random variable (r.v.), measurable function
  $ \ \zeta: \Omega \to R \ $ may be denoted alike

$$
||\zeta||_{q,\Omega} = ||\zeta||L(q,\Omega)  \stackrel{def}{=} \left[ \ \int_{\Omega} |\zeta(\omega)|^q  \ {\bf P}(d \omega) \right]^{1/q} = \left[ \ {\bf E}|\zeta|^q \ \right]^{1/q}, \ 1 \le q < \infty.
$$

\vspace{5mm}

 \ Let also  \ $ \  \xi_i(t) = \xi_i(t,\omega), \ t \in T,  \ \omega \in \Omega, \ i = 1, 2, \ldots, n, \ldots \ $ be a sequence of complete (total) measurable
 numerical valued {\it centered} (mean zero) $ \ {\bf E}\xi_i(t) = 0, \ t \in T \ $ separable  random fields (r.f.)
  or equally random processes (r.p.), defined in addition on our probability space $ \ (\Omega = \{\omega\}, B, {\bf  P}): \ $ \par

$$
\xi_i: T \otimes \Omega \to R.
$$

 \  Define

 \begin{equation} \label{Sn}
 S_n(t) \stackrel{def}{=} n^{-1/2} \sum_{i=1}^n \xi_i(t),
 \end{equation}

$$
R_n(t_1,t_2) := \Cov \left( S_n(t_1), S_n(t_2) \right) = n^{-1} \sum_{i=1}^n \sum_{j=1}^n {\bf E} \xi_i(t_1) \xi_j(t_2),
$$

$$
R(t_1,t_2) = R_{\infty}(t_1,t_2) := \lim_{n \to \infty}R_n(t_1,t_2),
$$
if there exists (and is finite) for all the values $ \ (t_1, t_2), \ t_1,t_2 \in T. \ $ \par
 \ Define the following r.f. $ \ S(t) = S(t,\omega) = S_{\infty}(t,\omega) \ $ as a Gaussian distributed separable centered r.f. with covariation function

\begin{equation} \label{limit pr}
{\bf E} S(t_1) S(t_2) = \Cov(S(t_1), S(t_2)) = R(t_1,t_2) = R_{\infty}(t_1,t_2).
\end{equation}

\vspace{5mm}

 \ {\bf Definition 1.1.} We will say as ordinary that the sequence of the r.f.  $ \ \{\xi_i(t) \} \ $  satisfies  the Central Limit Theorem (CLT) in the space
$ \  L(p,T) =L(p,T,\mu), \ $  iff all the r.f.  $ \ \{\xi_i(t) \} \ $  are (bi \ - ) measurable, belongs to the space $ \ L(p,T) \ $  with probability one:

$$
{\bf P} \left( ||\xi_i(\cdot)||L(p,T) < \infty  \right) = 1
$$
as well as the limited Gaussian r.f  $ \ S(\cdot), \ $  and  the sequence of normed r.f. $ \ S_n(\cdot) \ $ converges in distribution in this space $ \ L(p,T) \ $
 as $ \ n \to \infty \ $ to the r.f. $ \ S(\cdot): \ $ for every bounded continuous numerical valued functional $ \ G: L(p,T) \to R \ $

$$
\lim_{n \to \infty} {\bf E} G(S_n) = {\bf E} G(S).
$$

 \vspace{3mm}

  \ See the classical work of Yu.V.Prokhorov \cite{Prokhorov}. \par

 \ As a consequence:

\begin{equation} \label{norm distr}
\lim_{n \to \infty} {\bf P} ( ||S_n(\cdot)||L(p,T)| > u) = {\bf P} (||S(t)||L(p,T) > u), \ u > 0.
\end{equation}

\vspace{3mm}

 \ The asymptotic as well as non - asymptotic behavior for the right-hand side of the relation (\ref{norm distr}) as $ \ u \to \infty \ $ is
 investigated in particular in \cite{Dmitrovsky},  \cite{Piterbarg}. \par

\vspace{4mm}

 \ The CLT in different Banach spaces, or more generally, in linear topological ones,
  is devoted the extensive literature, see e.g. \cite{Acosta}, \cite{Araujo},
\cite{Dedecker}, \cite{Dudley},  \cite{ErDmOs}, \cite{ErOs}, \cite{Hoffmann}, \cite{KosOs},  \cite{Le Cam1}, \cite{Le Cam2},
\cite{Ledoux}, \cite{MalOs}, \cite{Mushtary},  \cite{Ostrovsky1}, \cite{Prokhorov}, \cite{Ratchckauskas}, \cite{Talagrand1}, \cite{Vakhania}. As a rule,  in  these works
was investigated the case when the r.f. $ \ \{\xi_i(t) \} \ $ are common  independent; the case when these fields are  weakly dependent
makes up the content of the articles \cite{ErDmOs}. \par

 \ The case of CLT in Lebesgue-Riesz spaces for weakly dependent variables under super strong mixing condition is considered in
\cite{ErOs}, \cite{MalOs}. \par

\vspace{4mm}

 \ The applications of the CLT in linear spaces in the statistics are described for instance, in \cite{Ostrovsky1}, \cite{Prokhorov}, \cite{Ratchckauskas}; in the theory of
 Monte-Carlo method-in  \cite{Frolov}, \cite{Grigorjeva}. \par

\vspace{4mm}

{\bf  We intend to generalize in this report  the foregoing results  \cite{ErOs}, \cite{MalOs} relative the CLT in the space $ \ L(2,T) \ $ into more general case of the space $ \ L(p,T) \ $
when the source random sequence $ \ \{\xi_i(t) \} \ $  forms relative the index $ \ i \ $ a sequence satisfying aside from certain moment assumptions  the strong mixing condition. } \par

\vspace{4mm}

 \ Recall some used definitions. Let $ \ B_1, \ B_2 \ $ be two sigma - subfields of source field $ \ B. \ $ The so - called {\it mixing coefficient} or equally Rosenblatt
 coefficient  $ \ \alpha(B_1, B_2) \ $   between $ \ B_1 \ $ and $ \ B_2 \ $ is defined as ordinary

$$
\alpha(B_1, B_2) \stackrel{def}{=} \sup_{A_1 \in B_1, \ A_2 \in B_2} | {\bf P}(A_1 \cap A_2) - {\bf P}(A_1) \ {\bf P}(A_2) |,
$$
see \cite{Rosenblatt}. \par

 \ Further, define a following family of sigma-algebras (fields)

$$
M_a^b \stackrel{def}{=} \sigma\{ \ \xi_i(t) \ \},   \ a \le i \le b, \ t \in T,
$$

\begin{equation} \label{alpha}
\alpha(i) \stackrel{def}{=}  \sup_n \max_{k \in [1,n]} \alpha( M_1^k, M_{k+i}^n).
\end{equation}

 \ By  definition, the sequence of r.f. $ \ \{\xi_i(t) \}, \ i = 1,2,\ldots \ $ satisfies a strong mixing condition, or equally $ \ \alpha \ $ mixing condition,
 iff $ \ \lim_{i \to \infty} \alpha(i)  = 0. \ $ \par
 Many examples of random processes and sequences obeys the strong mixing condition may be found in the articles \cite{Davydov}, \cite{Doukhan}. \par

 \vspace{3mm}

 \ There are many works devoted to the one-dimensional CLT  for such a sequences, including the invariance principle, see e.g. \cite{Hall}, \cite{Herrndorf},
 \cite{Ibragimov I},   \cite{Peligrad},  \cite{Rosanov}, \cite{Rosenblatt}, \cite{Utev1}, \cite{Utev2} etc.\par

 \ The non-asymptotical  estimation for sums of weakly dependent r.v. is represented in particular in \cite{ErDmOs}, \cite{ErOs}, \cite{Utev1}. We will apply
the estimate belonging to S.A.Utev  \cite{Utev2}. Let $ \ \{X_i\}, \ i = 1,2,\ldots,n \ $ be a sequence of centered one-dimensional r.v. satisfying the  mixing
condition with correspondent coefficient $ \ \alpha(i). \ $  Then when $ \ s = 2,4,6,\ldots \ $ and for every real value $ \ v \ge s \ $

\begin{equation} \label{Utev est}
{\bf E}|\sum_{i=1}^n X_i|^s  \le a_s \ \left[  \sum_{r=0}^{n-1} \alpha^{1 - s/v}(r) \ (r+1)^{s/2 - 1}  \right] \
\left( \ \sum_{i=1}^n {\bf E}^{2/v} |X_i|^v \ \right)^{s/2},
\end{equation}

where

$$
a_s = 12 \ (1 + 2s/3) \ (s-1) \ 3^s \ \frac{[s!]^2}{[(s/2)!]^2}.
$$

 \ It is no hard to evaluate

\begin{equation} \label{as estimate}
a_s^{1/s} \le K_U \cdot s, \ \hspace{4mm} K_U :=  2^{-5/12} \cdot 3 \cdot 7^{1/2} \cdot e^{ 2/e - 23/24} \approx 4.760327... ;
\end{equation}

\vspace{3mm}

 \ Let us denote

\begin{equation} \label{Z value}
Z[\alpha](s,v) := \left\{ a_s  \ \left[ \ \sum_{r=0}^{\infty} \alpha^{1 - s/v}(r) \ (r+1)^{s/2 - 1}  \ \right]  \ \right\}^{1/s},
\end{equation}
then if for some value $ \ v > s \ \Rightarrow   Z[\alpha](s,v) < \infty \ $  and $ \ X_i \in L(v,\Omega), \ $
we have under our conditions

\begin{equation} \label{Utev impr est}
||\sum_{i=1}^n X_i||_{s,\Omega}   \le Z[\alpha](s,v) \cdot \sqrt{ \sum_{i=1}^n ||X_i||_v^2}.
\end{equation}

 \ Define the variable

\begin{equation} \label{Y definition}
Y = Y[\{X_i\}](v) := \sup_i ||X_i||_{v,\Omega},
\end{equation}
then

\begin{equation} \label{via Y estimation}
\sup_n \ || \ n^{-1/2} \ \sum_{i=1}^n X_i||_{s,\Omega}  \le Z[\alpha](s,v) \ Y[\{X_i\}](v),
\end{equation}
 if of course the right-hand side is finite at last  for one value $ \ v; \ v > s. \ $ \par

\vspace{3mm}

\ If in particular all the centered variables $ \ \{X_i\} \ $ are in additional identical distributed, in particular, forms the
strictly stationary centered sequence, we conclude denoting $ \ X= X_1: \ $

\begin{equation} \label{Idv distr }
\sup_n \  ||n^{-1/2} \sum_{i=1}^n X_i||_{s,\Omega}   \le Z[\alpha](s,v) \cdot  ||X||_{v,\Omega}, \ v > s.
\end{equation}

\vspace{4mm}

 \ Of course,

\begin{equation} \label{Idv distr optim }
\sup_n \  ||n^{-1/2} \sum_{i=1}^n X_i||_{s,\Omega}   \le \inf_{v > s} \left\{ \ Z[\alpha](s,v) \cdot  ||X||_{v,\Omega} \ \right\}, \ v > s.
\end{equation}

\vspace{5mm}

\section{Moment estimates for the Lebesgue-Riesz norm of the sums of weakly dependent random fields.  }

\vspace{5mm}

 \ Let as before   $ \ \xi = \{ \xi_i(t) \}, \ t \in T, \ i = 1,2,\ldots \ $  be a sequence of  centered random fields, $ \ t \in T, \ $ satisfying the
 strong mixing condition  relative the index $ \ i \ $  with correspondent Rosenblatt coefficient (more exactly, the sequence of coefficients)
 $ \ \alpha = \{ \alpha(i) \ \}. \ $

 \ We get using the estimate (\ref{via Y estimation})

$$
\sup_n ||S_n(t)||_{s,\Omega} \le Z[\alpha](s,v) \ \sup_i ||\xi_i(t)||_{v,\Omega}, \ t \in T,
$$
or equally

$$
\sup_n {\bf E} |S_n(t)|^s \le  Z^s[\alpha](s,v) \ \left[ \ \sup_i ||\xi_i(t)||_{v,\Omega} \  \right]^s =
$$

$$
 Z^s[\alpha](s,v) \ \left[ \ \sup_i {\bf E}^{1/v} |\xi_i(t)|^v \ \right]^s, \ 1 \le s < v.
$$

\vspace{5mm}

 \ {\bf Proposition 2.1.}

\vspace{4mm}

 \ We deduce further after integration over the set $ \ T \ $ relative the measure $ \ \mu \ $ by virtue of theorem Fubini-Tonelli and Lyapunov-H\"older's inequality

\begin{equation} \label{key estim prelim}
\sup_n {\bf E} ||S_n||^s_{s,T} \le Z^s[\alpha](s,v) \ \int_T \left\{ \ \sup_i \left\{ \  {\bf E} |\xi_i(t)|^v \ \right\} \ \mu(dt) \ \right\}^{s/v}, \ 1 \le s < v.
\end{equation}
or equally

\begin{equation} \label{key estim}
\sup_n {\bf E} ||S_n||^s_{s,T} \le Z^s[\alpha](s,v) \  \sup_i {\bf E}^{s/v} ||\xi_i||^v_{v,T}, \ 1 \le s < v.
\end{equation}

\vspace{3mm}

 \ Of course, the last estimate is reasonable  if for instance the right-hand of (\ref{key estim}) is finite. \par

\vspace{5mm}

\section{Main result: Central Limit Theorem for the sums of weakly dependent random fields in the Lebesgue-Riesz space.  }

\vspace{5mm}

 \ {\bf Theorem 3.1.} Suppose that there exists a value $ \ v,   \ v > s, \ $ for which $ \  Z[\alpha](s,v) < \infty \ $ and
 $ \ \sup_i \left\{ \ \int_T {\bf E} |\xi_i(t)|^v \ \mu(dt)  \  \right\} < \infty. \ $ Suppose also that the sequence of r.v. $ \ \{\xi_i(\cdot)\} \ $
converges weakly (in distribution) as $ \ i \to \infty \ $ in the space $ \ L(s,\mu,T),  \ $ for instance, if all the r.v.  $ \ \{\xi_i(\cdot)\} \ $
are identical distributed.\par
 \ Then this sequence $ \ \{ \xi_i(t)\} \ $ satisfies CLT in the Lebesgue-Riesz space $ \ L(s,\mu,T). \ $ \par

 \vspace{5mm}

  \ {\bf Proof.} The convergence as $ \ n \to \infty \ $ of the characteristical functionals

$$
\phi_n(x) := {\bf E} \exp \left( i \int_T S_n(t) \ x(t) \ \mu(dt)   \right)
$$
 for arbitrary fixed non-random element  of conjugate space $ \ x(\cdot) \in L(s/(s-1), \mu,T) \ $ to the suitable for
the limiting Gaussian ones follows immediately from the one-dimensional limit theorems for strong random variables (Cramer's method), see e.g.
\cite{Peligrad},  \cite{Rosanov}, \cite{Rosenblatt}. The mentioned above one-dimensional r.v. haves the form

$$
\eta_i[x] =  \int_T \xi_i(t) \ x(t)  \ \mu(dt).
$$

\vspace{3mm}

 \ It remains to ground the weak compactness of distributions of the r.v. $ \ S_n(\cdot) \ $ in the space $ \ L(s,T). \  $ As long as the Banach space $ \ L(v,T,\mu) \ $
is separable,  and since the function $ \  y \to |y|^v \ $ satisfies the $ \ \Delta_2 \ $ condition,
there exists a compact linear operator $ \ U: L(v,T) \to L(v,T) \ $ such that the r.v. $ \  U^{-1} \xi_i(\cdot) \in L(v,T)  \  $ and moreover

$$
\sup_i \left\{ \ \int_T {\bf E} | U^{-1} [\xi_i](t)|^v \ \mu(dt)  \  \right\} < \infty,
$$
see  \cite{Ostrovsky5} and also \cite{Buldygin Supp}, \cite{Ostrovsky3}, \cite{Ostrovsky4}. One can apply  Proposition 2.1 for the sequence of r.f. $ \ U^{-1} [\xi_i](t): \ $

\vspace{4mm}

\begin{equation} \label{comp U estim}
\sup_n {\bf E} ||U^{-1} \ S_n||^s_{s,T} \le  W[\alpha](s,v),
\end{equation}
where

\begin{equation} \label{comp U estim via W}
W[\alpha](s,v) \stackrel{def}{=} Z^s[\alpha](s,v) \ \sup_i \left\{ \ \int_T {\bf E} |U^{-1} \ \xi_i(t)|^v \ \mu(dt)  \  \right\} < \infty, \ 1 \le s < v.
\end{equation}
 \ We conclude on the basis of Tchebychev's inequality for any value $ \ \epsilon \in (0,1) \ $

\begin{equation} \label{compactness}
\sup_n {\bf P} \left( \ ||U^{-1} S_n||L(s,T) > Y \ \right) \le \frac{W}{Y^s} \le \epsilon
\end{equation}
for sufficiently greatest positive value $ \ Y. \ $ \par

 \ Therefore, the sequence  of distributions of the r.v. $ \ S_n(\cdot) \ $ is  weak compact in the space $ \ L(s,T), \ $ \par
Q.E.D.\par

\vspace{4mm}

 \ {\bf Remark 3.1.} We have investigated above only the case when the number $  \ s \ $ is even:  $ \ s = 2,4,6,\ldots \ . \ $  Let now the number $  \ s \ $
be arbitrary greatest than 2: $ \ s > 2. \ $ Suppose in addition that the set $  \ T \ $ has a finite measure relative the measure $ \ \mu: \ \mu(T) = 1. \ $
Define the number $ \ \tilde{s} \ $ as a minimal even number  greatest than $ \ s: \ $

$$
\tilde{s} \stackrel{def}{=} \min\{  \ 2k, \ k = 1,2,\ldots; \ 2k \ge s \  \}.
$$

 \ As long as by virtue of Lyapunov's inequality $ \ ||f||L(s,T,\mu) \le ||f||(\tilde{s}, T,\mu), \ $ we conclude that if the conditions of
Theorem 3.1 are satisfied for the value $ \ \tilde{s} \ $ instead $ \ s, \ $ the CLT in the space $ \ L(s,T,\mu) \ $ holds true.\par

\vspace{5mm}

\section{ The case of superstrong mixing condition.}

\vspace{5mm}

 \ Recall that the so-called superstrong mixing coefficient between two sigma - algebras $ \ B_1,B_2 \subset B \ $ is defined by the formula

\begin{equation} \label{beta coeff}
\beta(B_1,B_2) \stackrel{def}{=} \sup_{{\bf P}(A_1 \in B_1) \ {\bf P}(A_2 \in B_2) > 0} \left| \ \frac{{\bf P}(A_1 \cap A_2) - {\bf P}(A_1) {\bf P}(A_2)}{ {\bf P}(A_1) {\bf P}(A_2)} \  \right|.
\end{equation}

 \  Define a following family of sigma-algebras (fields)

$$
M_a^b \stackrel{def}{=} \sigma\{ \ \xi_i(t) \ \},   \ a \le i \le b, \ t \in T,
$$
and consequently

\begin{equation} \label{alpha}
\beta(i) \stackrel{def}{=}  \sup_n \max_{k \in [1,n]} \beta( M_1^k, M_{k+i}^n).
\end{equation}

 \ By  definition, the sequence of r.f. $ \ \{\xi_i(t) \}, \ i = 1,2,\ldots \ $ satisfies a  super strong mixing condition,  or equally $ \ \beta \ $ mixing condition,
 iff $ \ \lim_{i \to \infty} \alpha(i)  = 0. \ $  This notion belongs to  Nachapetyan B.S. \cite{Nachapetyan} ; see also \cite{Bradley}. \par
 \ B.S.Nachapetyan in \cite{Nachapetyan}  proved that for the superstrong centered random sequence  $ \ \{X_i\}, \ i = 1,2,\ldots \ $

\begin{equation} \label{Nah est}
\sup_n || n^{-1/2} \sum_{i=1}^n X_i||L(s,\Omega) \le K_N [\beta](s) \ \sup_i ||X_i||L(s,\Omega), \ s \in [2,\infty),
\end{equation}
where

$$
K_N [\beta](s) := 2 s \left[ \ \sum_{k=1}^{\infty} \beta(k) \ (k+1)^{(s - 2)/2} \  \right]^{1/s}.
$$

\vspace{5mm}

 \ {\bf Theorem 4.1.} Suppose that the source sequence of centered random fields $ \ \{\xi_i(t)\} \ $ satisfies the superstrong mixing condition
 relative the index $ \ i, \ $

$$
\sup_i ||\xi_i(t)||L(s,\Omega) \in L(s, T,\mu), \ s \ge 2,
$$
 and

$$
K_N [\beta](s) < \infty.
$$

 \  Suppose also that the sequence of r.v. $ \ \{\xi_i(\cdot)\} \ $
converges weakly (in distribution) as $ \ i \to \infty \ $ in the space $ \ L(s,\mu,T),  \ $ for instance, if all the r.v.  $ \ \{\xi_i(\cdot)\} \ $
are identical distributed.\par

 \ Then this sequence $ \ \{ \xi_i(t)\} \ $ satisfies CLT in the Lebesgue-Riesz space $ \ L(s,\mu,T). \ $ \par

 \vspace{5mm}

  \ {\bf Proof} is much easier than one in theorem 3.1 and may be omitted. \par

\vspace{5mm}

\section{Concluding remarks.}

\vspace{5mm}

 \ {\bf A.} It is interest by our opinion to extend obtained in this report results into the the so-called anisotropic (Mixed) Lebesgue-Riesz spaces
 $ \ L(p_1,p_2, \ \ldots, p_d; \ T_1,T_2, \ldots, T_d), \ $ as well as to extend ones into the Sobolev's spaces. \par

 \ The independent case, i.e.when $ \ \alpha(r) = 0, \ r = 1,2, \ldots \ $ is investigated in \cite{Ostrovsky6}.\par

\vspace{4mm}

 {\bf B.} Perhaps, more general results in this direction  may be obtained by means of the so-called method of majorizing measures, see e.g.
\cite{Ledoux}, \cite{Talagrand1}. \par

\vspace{4mm}

{\bf C.} It is interest also by our opinion to deduce the non-asymptotical exponential decreasing  as  $ \ y \to \infty \ $ estimations for the tail of
distribution for the normed sums

$$
Q(y) := \sup_n {\bf P} (||S_n(\cdot)||L(s,T,\mu) > y), \ y \ge 1.
$$

\vspace{6mm}

\vspace{0.5cm} \emph{Acknowledgement.} {\footnotesize The first
author has been partially supported by the Gruppo Nazionale per
l'Analisi Matematica, la Probabilit\`a e le loro Applicazioni
(GNAMPA) of the Istituto Nazionale di Alta Matematica (INdAM) and by
Universit\`a degli Studi di Napoli Parthenope through the project
\lq\lq sostegno alla Ricerca individuale\rq\rq (triennio 2015-2017)}.\par

\end{document}